\documentclass[a4paper,12pt,reqno]{amsart}
\usepackage{amssymb}

\usepackage{amsmath}
\usepackage{amscd}
\usepackage{amsthm}
\usepackage{graphicx}
\usepackage[all]{xy}

\numberwithin{equation}{section}

\theoremstyle{plain}
 \newtheorem{thm}{Theorem}[section]
 \newtheorem{lem}[thm]{Lemma}

\theoremstyle{definition}
 \newtheorem{rem}{Remark}

 \newtheorem{defn}[thm]{Definition}

\begin{document}

\title{On the LMO conjecture}

\author{Takahito Kuriya}
\date{}
\address{Kyushu University, Faculty of Mathematics, 6-10-1 Hakozaki, Higasiku, Fukuoka 812-8581, Japan}
\email{marron@math.kyushu-u.ac.jp}
\maketitle

\begin{abstract} 
We give a proof of the LMO conjecture which say that for any simply connectd simple Lie group $G$, the LMO invariant of rational homology 3-spheres recovers the perturvative invariant $\tau^{PG}$. 
By Habiro-Le theorem, this implies that the LMO invariant is the universal quantum invariant of integral homology 3-spheres.
\end{abstract}

\section{Introduction}
In the late 1980s, E. Witten considered a quantum field theory whose Lagrangean is the Chern-Simons functional $CS(A)$, and proposed that the partition function of that provides a topological invariant of 3-manifolds. This path integral formula has not been justified mathematically. However, there were two approaches, called the operator formalism and the perturbative expansion, to obtain observables from the path integral. On one hand, in a mathematically rigorous way based on the operator formalism, Reshetikhin and Turaev constructed invariants of 3-manifolds, which are now called quantum invariants in \cite{RT}. T. Ohtsuki considered the perturbative invariant $\tau^{SO(3)}$ of 3-manifolds obtained from the quantum invariant $\tau^{SO(3)}_r$ by taking $r$-adic limit (general cases are due to T.T.Q. Le). On the other hand, as a mathematical counterpart of the perturbative expansion of the Chern-Simons path integral, T.T.Q. Le, J. Murakami, and T. Ohtsuki {\cite{LMO}} constructed an invariant, called the LMO invariant, of 3-manifolds, and gave a conjecture (the LMO conjecture) which says that for any simply connectd simple Lie group $G$, the LMO invariant of rational homology 3-spheres recovers the perturvative invariant $\tau^{PG}$ (see Remark 7.5 in {\cite{LMO}}). It is very significant to prove this conjecture from the viewpoint that it is related to the problem whether we obtain the same observables from the path integral by the above mentioned two approaches. T.T.Q. Le and K. Habiro's theorem implies that the LMO invariant is the universal quantum invariant of integral homology 3-spheres. 

When $G=SU(2)$, the LMO conjecture has been shown by T. Ohtsuki {\cite{Oh2}}. The proof of the main result differs from that of Ohtsuki for the case $G=SU(2)$, since T. Ohtsuki \cite{Oh2} used some identities which are either specific to the $G=SU(2)$ case or hard to generalize to the other cases. It is announced in {\cite{BGRT}} (Corollary 2.14) that the LMO conjecture has been proved, but the proof is not published yet. 

The goal of this paper is to prove the LMO conjecture. Our proof is based on a general formula of the perturvative invariant $\tau^{PG}$ and the $\r{A}rhus \ integral$. The former (Theorem \ref{PG}) was given by T.T.Q. Le in {\cite{Le2}} and the latter (see Section \ref{s-LMO}) was defined by D. Bar-Natan, S. Garoufalidis, L. Rozansky and D.P. Thurston in {\cite{BGRT}}. The $\r{A}rhus \ integral$ gives a reconstruction of the LMO invariant. However, at least the auther, it seems to be easier to define and to compute. We should note that the definition of the LMO invariant works on arbitrary closed oriented 3-manifolds but the $\r{A}rhus \ integral$ does only on rational homology 3-spheres, although this is sufficient for our purpose. In addition, we use $\mathfrak{g}$-weight system translates graphical derivative into the usual one. In \cite{TK}, the auther proved that if the LMO conjecture for an arbitrary rational homology 3-sphere which can be obtained by integral surgery along some algebraically split framed link is true, then the LMO conjecture is true (Theorem \ref{37}). Using these facts, we will prove the following theorem.

\begin{thm}[LMO conjecture, BGRT or Le theorem]{\label{27}}
For any simply connected compact simple Lie group $G$, 
$\hat{Z}^{\mathrm{LMO}}$ recovers $\tau^{PG}$.
\end{thm}

The paper is organized as follows. In Section \ref{s-LMO} we review some notations related to the LMO invariant. In Section \ref{WS} we recall the definition of the $\mathfrak{g}$-weight system. In Section \ref{s-Q} we recall the Le's theorem (Theorem \ref{PG}). Section \ref{p-LMO} is devoted to proof of the LMO conjecture.\\

{\bf Acknowledgment.}
I would like to thank many people for their help with this paper, especially, to Professor M. Wakimoto, Professor K. Habiro and Professor T. Ohtsuki for their helpful comments and to Professor Mitsuyoshi Kato and Professor Osamu Saeki for their encouragement.

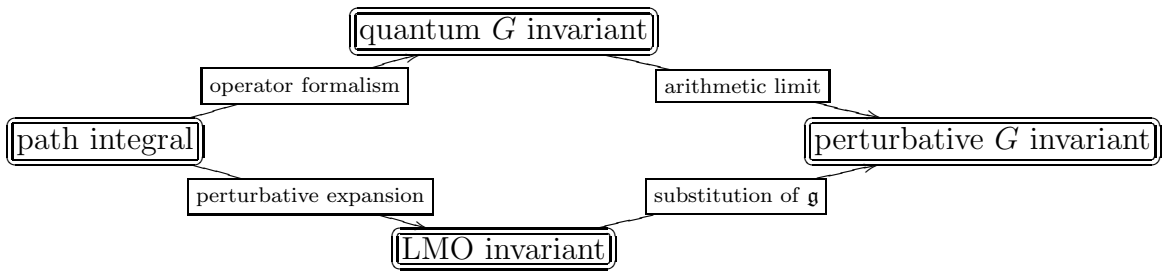
\begin{figure}[ht]
\[ 
\xymatrix{ 
 & & *+[F=:<3pt>]{\text{quantum $G$ invariant}} \ar[drr]|-{\small{\text{\fbox{arithmetic limit}}}} & \\
 *+[F=:<3pt>]{\text{path integral}} \ar[urr]|-{\small{\text{\fbox{operator formalism}}}} \ar[drr]|-{\small{\text{\fbox{perturbative expansion}}}} & & & & *+[F=:<3pt>]{\text{perturbative $G$ invariant}} \\
& & *+[F=:<3pt>]{\text{LMO invariant}} \ar[urr]|-{\small{\text{\fbox{substitution of $\mathfrak{g}$}}}} 
} 
\] 
\caption{Related invariants of rational homology 3-spheres in this paper}
\end{figure}

\section{The LMO invariant}{\label {s-LMO}}

A framed link is the image of an embedding of a disjoint union of annuli into $\mathbb{R}^3$ or $\mathbb{S}^3$. The underlying link of a framed link is the link obtained by restricting an annulus $\mathbb{S}^1 \times [0,1]$ to its center line $\mathbb{S}^1 \times \{\frac{1}{2}\} $. The framing of a component of a framed link is the isotopy class of framed knots whose underlying knots are equal to the component.\\
 A uni-trivarent graph is a graph every vertices of which is either univarent or trivarent. A uni-trivarent graph is vertex-oriented if at each trivarent vertex a cyclic order of edges is fixed. A 3-valent (resp.1-varent) vertex is called an internal (external) vertex. Let $X$ be a compact oriented 1-dimensional manifold. A \textit{chord diagram} with support $X$ is the manifold $X$ together with a vertex-oriented uni-trivalent graph whose external vertices are on $X$;and the graph dose not have any connected component homeomorphic to a circle. In figures components of $X$ are depicted by solid lines, while the graph is depicted by dashed lines, with the convention that the orientation at every vertex is counterclockwise. There may be connected components of the dashed graph which do not have univalent vertices, ant hence do not connect any solid lines.\\
 All vector spaces are over the field $\mathbb{Q}$ of rational numbers. Let $\mathcal{A}(X)$ be the vector space spanned by chord diagrams with support $X$, subject to the $AS$,$IHX$ and $STU$ relation.

\begin{figure}[ht]
\begin{center}
\begin{minipage}{1.5cm}
\includegraphics[width=\hsize]{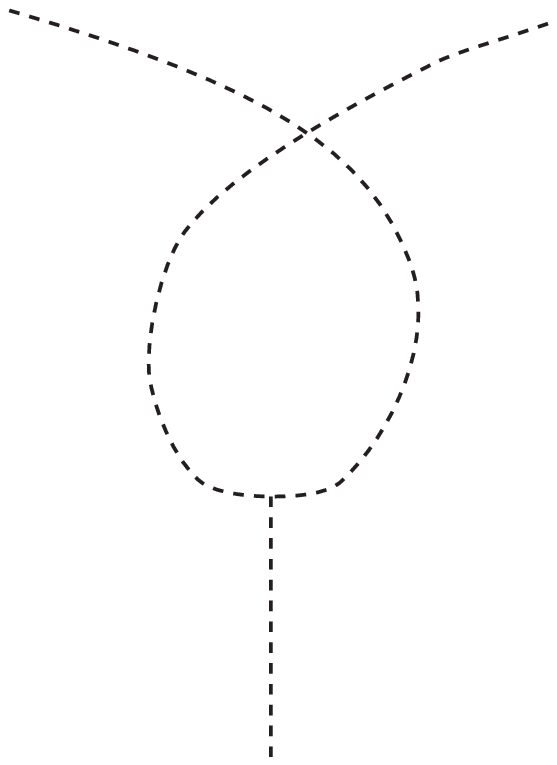}
\end{minipage}
$=\ -$
\begin{minipage}{1.5cm}
\includegraphics[width=\hsize]{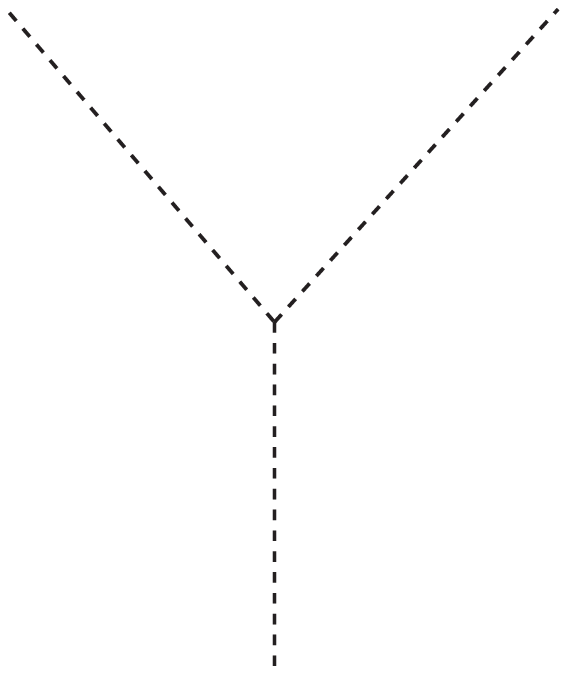}
\end{minipage}
\end{center}
\caption{AS relation:}
\end{figure}

\begin{figure}[ht]
\begin{center}
\begin{minipage}{1.3cm}
\includegraphics[width=\hsize]{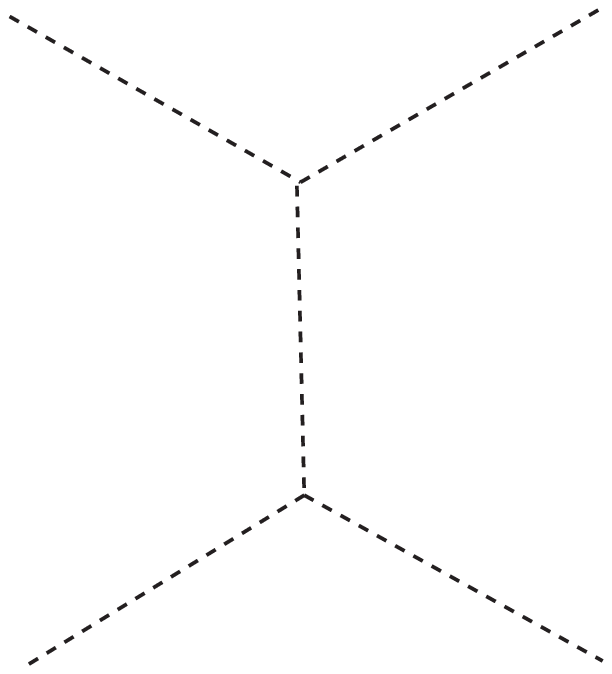}
\end{minipage}
$=$
\begin{minipage}{1.5cm}
\includegraphics[width=\hsize]{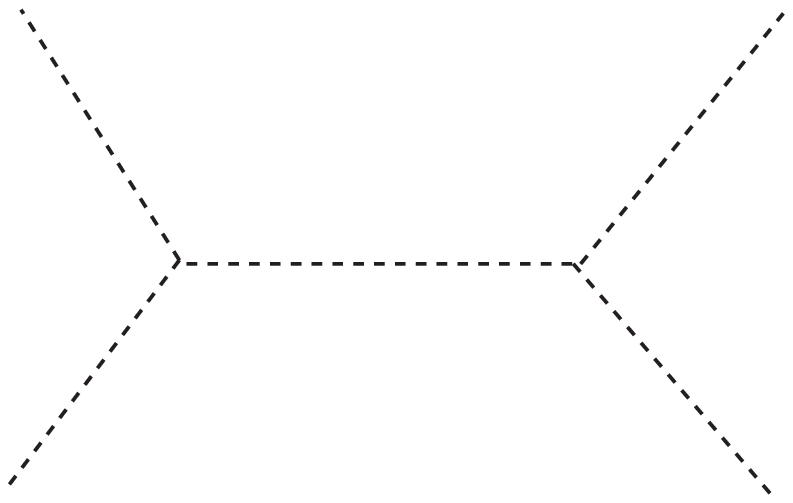}
\end{minipage}
$-$
\begin{minipage}{1.5cm}
\includegraphics[width=\hsize]{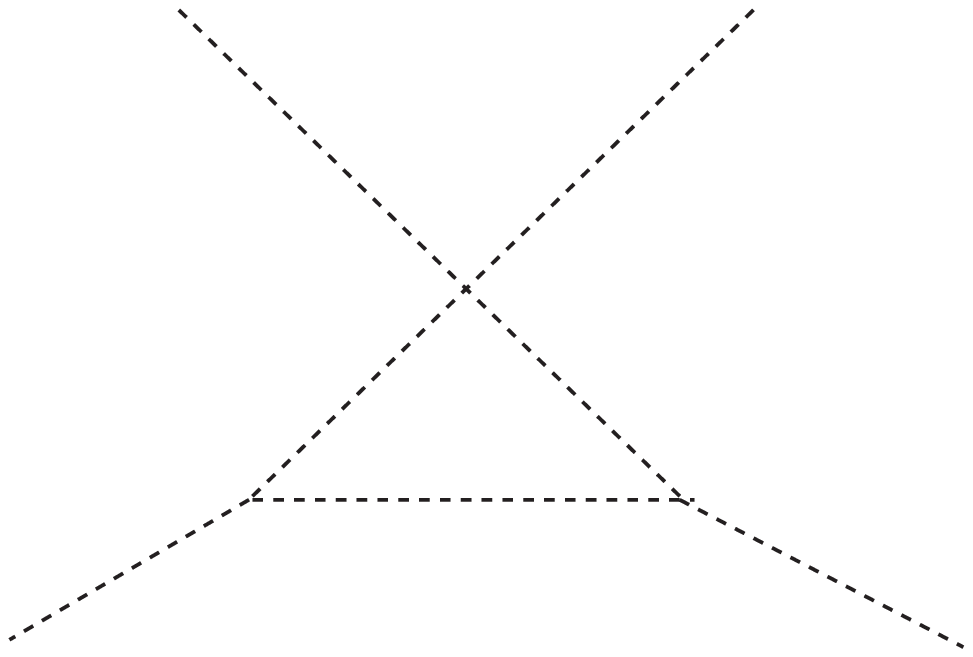}
\end{minipage}
\end{center}
\caption{IHX relation:}
\end{figure}

\begin{figure}[ht]
\begin{center}
\begin{minipage}{1.5cm}
\includegraphics[width=\hsize]{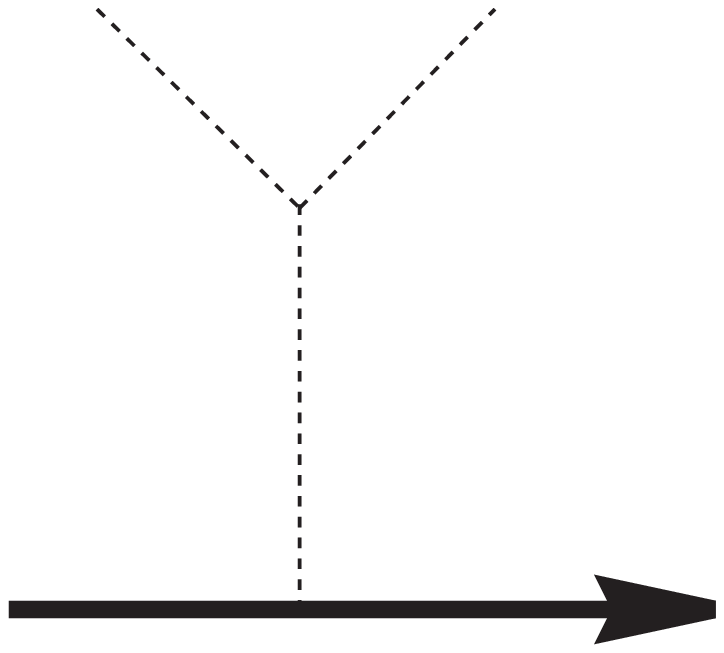}
\end{minipage}
$=$
\begin{minipage}{1.5cm}
\includegraphics[width=\hsize]{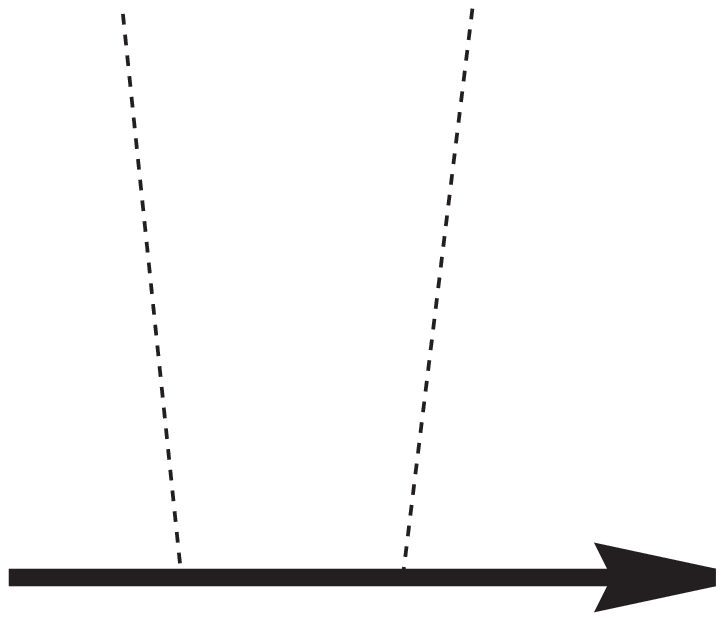}
\end{minipage}
$-$
\begin{minipage}{1.5cm}
\includegraphics[width=\hsize]{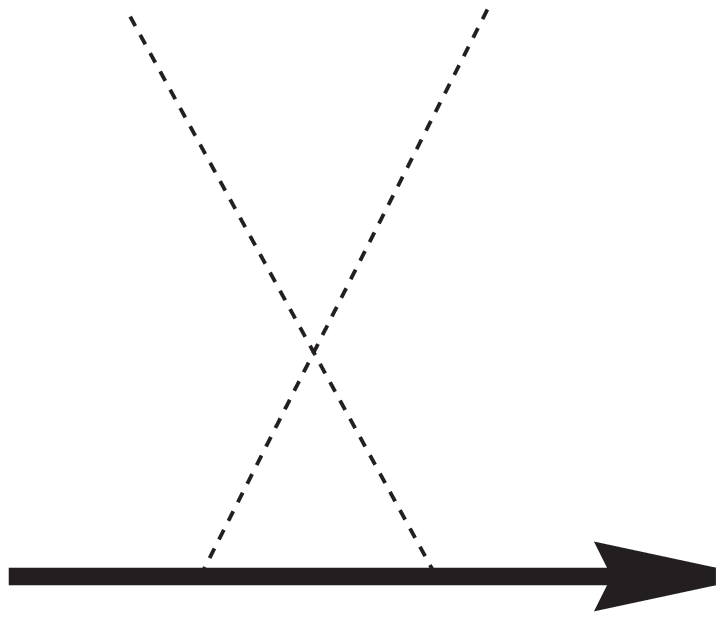}
\end{minipage}
\end{center}
\caption{STU relation}
\end{figure}
 The degree of a chord diagram is half the number of vertices of the dashed graph. We also use $\mathcal{A}(X)$ to denote the completion of $\mathcal{A}(X)$ with respect to the grading.\
 Of particular interest is the space $\mathcal{A}(\phi)$, i.e. when there are no solid lines nor circles. Then every dashed graph must be 3-valent. This vector space $\mathcal{A}(\phi)$ is a commutative algebra in which the product of 2 graphs is just the union of them.\
 A general way to remove solid loops. Let $\mathcal{A}(m)$, for any positive integer m, be the vector space spanned by uni-trivarent vertex-oriented graphs with exactly m 1-valent vertices numbered 1,2,...,m, subject to the $AS$ and $IHX$ relations.\
 $\mathcal{A}(\mathbb{S}^1)$ is isomorphic to $\mathcal{B}$, the space of uni-trivalent graphs (subject to the $IHX$,$AS$ conditions). This isomorphism $\sigma :\mathcal{A}(\mathbb{S}^1) \to \mathcal{B}$ is the formal PBW linear isomorphism between $\mathcal{A}(\mathbb{S}^1)$ and $\mathcal{B}$.The map $\sigma$ is most easily defined as the inverse of the symmetrization map $\chi:\mathcal{B} \to \mathcal{A}(\mathbb{S}^1)$. Formally, $\displaystyle{\mathcal{B} \ =\ \prod ^{\infty }_{m=0}\mathcal{A}(m)/S_m}$, where the symmetric group $S_m$ act on $\mathcal{A}(m)$ by permuting the external vertices. $\mathcal{B}$ is an algebra where the product is the disjoint union $\sqcup $.

\begin{defn}{\label{4}}
Let $\Omega \in \mathcal{B}$ be given by\\
\begin{center}
$\displaystyle{\Omega=\textrm{exp}_{\sqcup}\sum^{\infty }_{m=1}b_{2m}\omega_{2m}}$
\end{center}
Here the constants $b_{2m}$ are the modified Bernoulli numbers, defined by the power series expansion \\
\begin{center}
$\displaystyle{\sum^{\infty }_{m=0}b_{2m}x^{2m}=\frac{1}{2}\textrm{log}\frac{\textrm{sinh}\frac{1}{2} x}{\frac{x}{2}}}$
\end{center}
and $\omega_{2m}$ denotes the 2m-wheel,the degree 2m uni-trivalent graph made of a 2m-gon with 2m legs.
\end{defn}

\begin{defn}{\label{5}}
For any element $D \in \mathcal{B}$, define a map $\hat{D}=\partial D:\mathcal{B} \to \mathcal{B}$ \\
to act on diagrams $D' \in \mathcal{B}$ by gluing all legs of $D$ to some subset of legs of $D'$. Likewise, define a pairing 
\begin{center}
\[
<D,D'>:=
\left\{
   \begin{array}{ll}
sum \ of \ all \ ways \ of \ gluing \\ the \ legs \ of \ D \ to \ the \ legs \ of \ D' & \mbox{(if $\sharp$ (legs \ of \ D) \ =\ $\sharp$ (legs \ of \ D')\ )} \\
   0 & \mbox{(otherwise)}
   \end{array}
\right.
\]\\
\end{center}

\end{defn}

\begin{defn}{\label{6}}

Suppose that $K$ is a Morse knot, that is a 
knot embedded in $\mathbb{R}^3=\mathbb{R}_t\times \mathbb{C}_z$ such that the projection to the $\mathbb{R}$
coordinate is a Morse function. 

Then we define the element in $\mathcal{A}(\mathbb{S}^1)/<\Theta>$ by
$$\overline z (K):=\sum_{m=0}^\infty \frac{1}{(2\pi i)^m}
\int\limits_{\genfrac{}{}{0pt}{1}{t_{\text{min}}<  t_1 < \dots <
t_m< t_{\text{max}}}{t_j \text{ non-critical}}}
\sum_{\genfrac{}{}{0pt}{1}{\text{pairings }}{P=(\{z_j,z_j'\})}}(-1)^{\downarrow P}
D_P\bigwedge_{j=1}^m \frac{dz_j-dz'_j}{z_j-z'_j}$$

where a set of pairings $P$ is a set of unordered pairs
$\{z_j,z_j'\}$ so that 
$(t_j,z_j),(t_j,z_j')\in \mathbb{R} \times \mathbb{C}$ are distinct points on the knot
on the same level,  where $\!\downarrow\! P$ is the number of points
in the pairing such that the orientation is pointing downwards and
where $D_P$ is the diagram obtained by drawing the preimage of the
knot and drawing chords connecting the points corresponding to the pairs. $\Theta$ denotes the $\lq\lq$isolated chord" diagram in $\mathcal{A}(\mathbb{S}^1)$.
Considered as an element of $\mathcal{A}(\mathbb{S}^1)$, $\overline z (K)$ gives an invariant of Morse
knots, but is not invariant if one 
introduces a minima-maxima pair with respect to the
$\mathbb{R}$-axis, which increases the number of critical points.  
It is necessary to introduce the following correction
term:  define $\nu:=\overline
z\left( U \right)^{-1}\in \mathcal{A}(\mathbb{S}^1)$,
where U is a Morse embedding of the unknot with four critical
points.  Now the
Kontsevich integral for a Morse knot $K$ with $c(K)$ critical points
is defined to be
$$Z(K):=\nu^{c(K)/2}\overline z(K) \in \mathcal{A}(\mathbb{S}^1).$$

\end{defn}

\begin{defn}{\label{7}}
The integer framed version of Kontsevich integral is defined as follows:\\ 
Let $K$ be a integer framed knot, with  framing $f$.  Set
\begin{center}
$\displaystyle{\hat{Z}(K):=Z(K)\textrm{exp}(\frac{f}{2}\Theta ))}$
\end{center} 
where in this equation, $\Theta$ denotes the $\lq\lq$isolated chord" diagram in $\mathcal{A}(\mathbb{S}^1)$.

\begin{figure}[ht] 
\begin{center}
$\Theta =$
\begin{minipage}{1.5cm}
\includegraphics[width=\hsize]{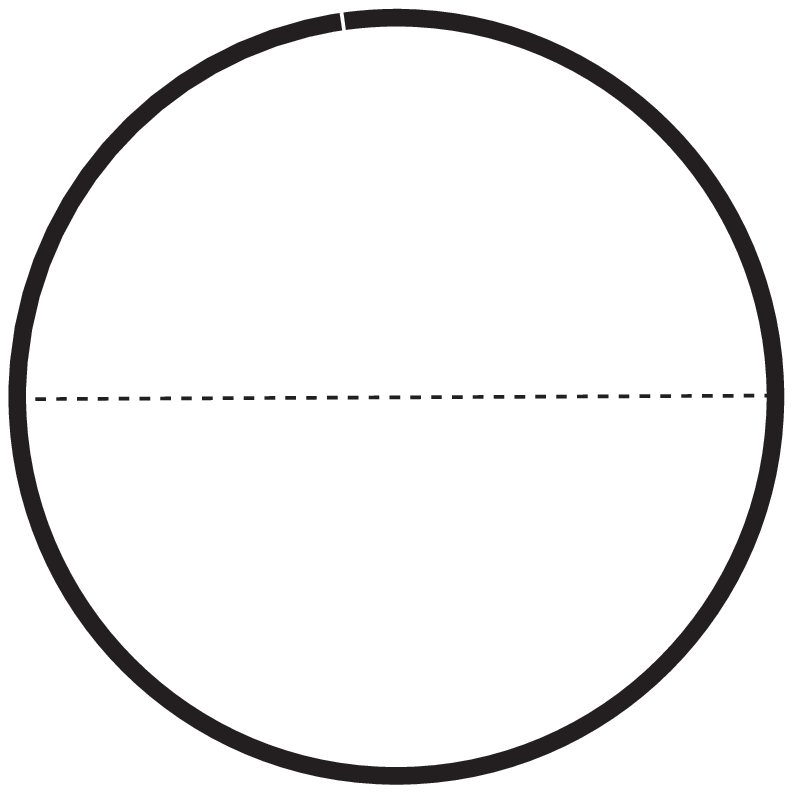}
\end{minipage}
\end{center}
\caption{$\Theta$ graph}
\end{figure}
\end{defn}

\begin{defn}{\label{10}}
The LMMO version of the Kontsevich integral of framed knot $\check{Z}$ is defined by 
\begin{center}
$\displaystyle{\check{Z}(K)= \nu \hat{Z}(K)}$.
\end{center}
\end{defn}

Since the framing is $f$, one has that (in $\mathcal{B}$) 
\begin{center}
$\displaystyle{\sigma \check{Z}(K)=\textrm{exp}_{\sqcup}(\frac{f}{2} \omega_1) \sqcup Y}$
\end{center}
where $\omega_1$ is the $\lq\lq$dashed interval" (without internal vertices), and $Y$ is an element in $\mathcal{B}$ every term of which must have at least one internal vertex.

\begin{defn}{\label{11}}
The formal Gaussian integral of the Kontsevich integral is defined by
\begin{center}
$\displaystyle{\int^{FG} \sigma \check{Z}(K)=\langle\textrm{exp}_{\sqcup}(-\frac{1}{2f} \omega_1),Y\rangle}$.
\end{center}
\end{defn}

For simplicity we assume that $M$ is obtained from $\mathbb{S}^3$ by surgery along a framed knot $K$. The case of links is similar (see {\cite{BGRT}}).

\begin{defn}{\label{12}}
The LMO invariant $\hat{Z}^{\mathrm{LMO}}(M)$ of a rational homology 3-sphere $M$ is defined by
\begin{center}
$\displaystyle{\hat{Z}^{\mathrm{LMO}}(M):=\frac{\int^{FG} \sigma \check{Z}(K)}{\int^{FG} \sigma \check{Z}(\bigcirc ^{\text{sign}(f)})}}$
\end{center}
\end{defn}

\begin{lem}{\label{LMOFG}}
Let $M_K$ be a rational homology 3-sphere obtained from $\mathbb{S}^3$ by surgery along a framed knot $K$ with framing $f$ and let $K_0$ be the same knot with framing $0$. Then
$$\hat{Z}^{\mathrm{LMO}}(M_K)=\langle\Omega,\Omega\rangle\exp \big( \frac{3\mathrm{sign}(f)-f}{48}\theta \big)\int^{FG} \big(\hat{\Omega}^{-1}\sigma \check{Z}(K_0)\big)\exp\frac{f}{2}\omega_1$$
where in this equation, $\theta$ denotes the trivalent graph in $\mathcal{A}(\phi)$.
\end{lem}

\begin{proof}
In \cite{{BGRT},{BL}}, the following facts were proved.
\begin{xalignat*}{2}
\int^{FG} \hat{\Omega}^{-1}D=&\int^{FG}D & \\
\int^{FG} \sigma \check{Z}(\bigcirc ^{\text{sign}(f)})=&\langle\Omega,\Omega\rangle^{-1}\exp \big(-\frac{\mathrm{sign}(f)}{16}\theta \big) & \\*
\hat{\Omega}^{-1}\sigma \Theta=&\ \omega_1-\frac{1}{24}\theta &
\end{xalignat*}
Note $\hat{\Omega}^{-1}\sigma$ is an algebra isomorphism.
\end{proof}

\section{weight system}{\label {WS}}
Let us recall the definition of the $\mathfrak{g}$-weight system. The map $\hat{W}_{\mathfrak{g}}$ is a combination of the two maps below. 
First, a graph $D\in \mathcal{A}(\mathbb{S}^1)$ (or $\mathcal{B}$) is multiplied by $h^{\mathrm{deg}(D)}$. Second map is as follows.
Let $\mathfrak{g}$ be a Lie algebra with an $\mathfrak{g}$-invariant inner product
$b: \mathfrak{g} \otimes \mathfrak{g} \to \mathbb{C}$. To each chord diagram $D$ with $m$ univalent
vertices  we assign a tensor 
$$
   T_{\mathfrak{g},b}(D) \in \mathfrak{g}^{\otimes m}
$$ 
as follows.

The Lie bracket $[\ ,\ ]: \mathfrak{g} \otimes \mathfrak{g} \to \mathfrak{g}$ can be considered
as a tensor in $\mathfrak{g}^*\otimes \mathfrak{g}^* \otimes \mathfrak{g}$. The inner product 
$b$ allows us to identify the $\mathfrak{g}$-modules 
$\mathfrak{g}$ and $\mathfrak{g}^*$, and therefore  $[\ ,\ ]$ can be considered as a tensor
$f\in (\mathfrak{g}^*)^{\otimes 3}$ and $b$ gives rise to an invariant symmetric
tensor $c \in \mathfrak{g} \otimes \mathfrak{g}$.

For a chord diagram $D$ denote by $T$ the set of its trivalent
vertices, by $U$ the set of its univalent (exterior) vertices, and by
$E$ the set of its edges. Taking $|T|$ copies of the tensor $f$ and $|E|$
copies of the tensor $c$ we consider a new tensor
$$
\widetilde{T}_\mathfrak{g}(D) = \Bigl(\bigotimes_{v\in T} f_v\Bigr) \otimes
                     \Bigl(\bigotimes_{\ell\in E} c_\ell\Bigr)
$$
which is an element of the tensor product
$$
\mathfrak{g}^D = 
\Bigl( \bigotimes_{v \in T}
(\mathfrak{g}^*_{v,1}\otimes \mathfrak{g}^*_{v,2} \otimes \mathfrak{g}^*_{v,3})
\Bigr)
\otimes
\Bigl( \bigotimes_{\ell \in E}
(\mathfrak{g}_{\ell,1}\otimes \mathfrak{g}_{\ell,2})
\Bigr),
$$
where $(v,i), \ i=1,2,3,$ mark the three edges meeting at the
vertex $v$ (consistently with the cyclic ordering of these edges), and
$(\ell,j), \ j=1,2,$ denote the endpoints of the edge $\ell$. 

Since $c$ is symmetric and $f$ is completely antisymmetric,
the tensor $\widetilde{T}_\mathfrak{g}(D)$ does not depend on the choices of
orderings.

If $(v,i)=\ell$ and $(\ell,j)=v$, there is a natural contraction map
$\mathfrak{g}^*_{v,i} \otimes \mathfrak{g}_{\ell,j} \to \mathbb{C}$. Composition of all such 
contractions gives us a map
$$
\mathfrak{g}^D \longrightarrow  \bigotimes_{u \in U} \mathfrak{g} = \mathfrak{g}^{\otimes m}, \ \text{where} \ m=|U|.
$$

The image of $\widetilde{T}_\mathfrak{g}(D)$ in
$\mathfrak{g}^{\otimes m}$ is denoted by $T_{\mathfrak{g},b}(D)$ (or usually just by
$T_{\mathfrak{g}}(D)$).

If $D\in \mathcal{A}(\mathbb{S}^1)$, then the STU relations project $T_{\mathfrak{g}}(D)$ to $W_{\mathfrak{g}}(D) \in \mathcal{U}(\mathfrak{g})^{\mathfrak{g}}$. If $D\in \mathcal{B}$, then the equivalences between univalent vertices project $T_{\mathfrak{g}}(D)$ to $W_{\mathfrak{g}}(D) \in \mathcal{S}(\mathfrak{g})^{\mathfrak{g}}$. The conditions STU and AS,IHX are automatically
satisfied for $W_\mathfrak{g}$: the relations AS,IHX are 
the anticommutativity 
and the Jacobi identity for the Lie bracket, and STU is
just the definition of the universal enveloping algebra as a quotient
of the tensor algebra of $\mathfrak{g}$.

$$
W_\mathfrak{g}: \mathcal{A}(\mathbb{S}^1) \to Z(\mathcal{U}(\mathfrak{g})) 
$$
is called {\em the universal weight system}
corresponding to $\mathfrak{g}$. It is universal in the sense that any $W_{\mathfrak{g},V}$
constructed using a representation $V$ of the Lie algebra $\mathfrak{g}$ (see
{\cite{BN}})  is an evaluation of $W_\mathfrak{g}$:
$$
W_{\mathfrak{g},V}(D) = \text{Tr}_{V}\bigl(W_\mathfrak{g}(D)\bigr).
$$

\section{the quantum invariant and the perturbative invariant}{\label {s-Q}}

Let $G$ be a simply connected Lie group, $\mathfrak{g}$ its Lie algebra, and $V_\lambda$ finite dimensional irreducible representations of $\mathfrak{g}$ (and those of $U_q (\mathfrak{g})$ for generic $q$) parametrized  by dominant integral weights $\lambda$. Given a simple complex Lie algebra $\mathfrak{g}$ with an invariant inner product $(\ ,\ )$, As we mentioned before, there is a map 
$$\hat{W}_{\mathfrak{g}}:\mathcal{B} \to \mathcal{S}(\mathfrak{g})^{\mathfrak{g}}[[h]]$$
We call the map 
$$\hat{W}_{\mathfrak{g}}\circ \sigma \circ \hat{Z}: \{ \text{Knots} \}\to \mathcal{S}(\mathfrak{g})^{\mathfrak{g}}[[h]]$$
the $\mathfrak{g}$-colored Jones function of a knot. The $\mathfrak{g}$-colored Jones function of a knot is a generating function of the quantum invariants of a knot. Given an irreducible representation $V_{\lambda}$ of $\mathfrak{g}$, its evaluationat its dominant weight gives rise to a linear map 
$$\mathcal{S}(\mathfrak{g})^{\mathfrak{g}}[[h]] \to \mathbb{Q}[[h]]$$
As explained in {\cite{Ro2}}, after multipling $\mathrm{dim}V_{\lambda}$, the image of the Kontsevich integral under the composition above two maps is an element of the ring $\mathbb{Z}[q^{\pm 1}]$ (where $q=e^h$) and coincides with the quantum group invariant $J_{\mathfrak{g},V_{\lambda}}(K)$ of knot $K$, using the $(\mathfrak{g},V_{\lambda})$ deta. We denote $J_{\mathfrak{g},V_{\lambda}}(K)\times J_{\mathfrak{g},V_{\lambda}}(\text{unknot})$ by $Q_\mathfrak{g}(K)(\lambda)$. From $Q_\mathfrak{g}(K)$, one can construct the projective quantum invariant of rational homology 3-spheres , we denote by $\tau^{PG}_r(M)$. $\tau^{PG}_r(M)$ belongs to $\mathbb{C}$, which is not polynomial ring nor a power series ring. We can not expand $\tau^{PG}_r(M)$ with respect to $r$ in an ordinary sense. Insted of expansion in $r$, T. Ohtsuki considered the perturbative invariant $\tau^{SO(3)}$ of 3-manifolds obtained from the quantum invariant $\tau^{SO(3)}_r$ by taking $r$-adic limit (general cases were defined by T.T.Q. Le, see {\cite{Le2}}).

In \cite{Le2}, T.T.Q. Le showed below. 
\begin{thm}{\label{PG}}
Let $M_K$ be a rational homology 3-sphere obtained from $\mathbb{S}^3$ by surgery along a framed knot $K$ with framing $f$ and let $K_0$ be the same knot with framing $0$. Then
$$\tau^{PG}(M_K)=\frac{1}{|W|}q^{\frac{\mathrm{sign}(f)-f}{2}|\rho|^2}\prod_{\alpha>0} (1-q^{\mathrm{sign}(f)(\rho,\alpha)}) \sum_{\beta \in Y ,n\in \mathbb{Z}_+ \atop {2|\Phi_+|\le 2j\le n+2|\Phi_+|}} c_{\beta,2j,n}(2j-1)!! \big(-\frac{|\beta|^2}{f} \big)^j h^{n-j} $$
where $\displaystyle{\big( Q_\mathfrak{g}(K_0)\big|_{q=e^h} \big)(\lambda-\rho)=\sum_{\beta \in Y ,n\in \mathbb{Z}_+ \atop {2|\Phi_+|\le j\le n+2|\Phi_+|}} c_{\beta,j,n}\beta^j(\lambda)h^n}$ , $\beta^j(\lambda)=(\beta,\lambda)^j$ , $Y$ is the root lattice , $\Phi_+$ is the set of all positive roots and $\rho=\frac{1}{2}\sum_{\alpha >0}\alpha $.
\end{thm}

{\section {Proof of the LMO conjecture}}{\label{p-LMO}}

In \cite{TK}, the auther proved the next theorem.

\begin{thm}{\label{37}}
If the LMO conjecture for an arbitrary rational homology 3-sphere which can be obtained by integral surgery along some algebraically split framed link is true, then the LMO conjecture is true.
\end{thm}

In an old version (\cite{G0}) of \cite{G1}, S. Garoufalidis proved the following fact.

\begin{thm}{\label{CM}}
The following diagram commutes
\[
\begin{CD}
 \mathcal{B} @>{\ \ \ \ \hat{\Omega}\ \ \ \ \ \ \ }>> \mathcal{B}\\
@V{\hat{W}_{\mathfrak{g}}}VV         @V{\hat{W}_{\mathfrak{g}}}VV  \\
\mathcal{S}(\mathfrak{g})^{\mathfrak{g}}[[h]] @>>{\ \ \ \ D(j_{\mathfrak{g}} ^{\frac{1}{2}})\ \ \ \ }>  \mathcal{S}(\mathfrak{g})^{\mathfrak{g}}[[h]]
\end{CD} 
\]
where $D(j_{\mathfrak{g}} ^{\frac{1}{2}})$ is the Duflo isomorphism of commutative algebras. 
\end{thm}

\begin{rem}{\label{Du}}
The only fact that we should note is that $D(j_{\mathfrak{g}} ^{\frac{1}{2}})$ is shift by $\rho$.
\end{rem}

\begin{proof}[Proof of Theorem \ref{27}]
In \cite{G0}, S. Garoufalidis proved the following formula.
\begin{xalignat}{2}{\label{Omega}}
\hat{W}_{\mathfrak{g}}(\langle\Omega,\Omega\rangle)=&\prod _{\alpha >0} \frac{\text{sinh} \frac{(\rho,\alpha)}{2}h} {\frac{(\rho ,\alpha )}{2}h} & 
\end{xalignat}
For any simply connected compact simple Lie group $G$, we have
\begin{xalignat}{2}{\label{Theta}}
W_\mathfrak{g} \big( \theta \big)=&24(\rho,\rho) & {\text{(see \cite{TK})}} 
\end{xalignat}

\begin{xalignat*}{2}
&\hat{W}_{\mathfrak{g}}\big( \hat{Z}^{\mathrm{LMO}}(M_K) \big) & \\
=&\hat{W}_{\mathfrak{g}}\big(\langle\Omega,\Omega\rangle\exp \big( \frac{3\mathrm{sign}(f)-f}{48}\theta \big)\int^{FG} \big(\hat{\Omega}^{-1}\sigma \check{Z}(K_0)\big)\exp\frac{f}{2}\omega_1\big) & \\
=&\hat{W}_{\mathfrak{g}}(\langle\Omega,\Omega\rangle)\exp \big( \frac{3\mathrm{sign}(f)-f}{48}\hat{W}_{\mathfrak{g}}(\theta) \big) \hat{W}_{\mathfrak{g}}\Big(\int^{FG} \big(\hat{\Omega}^{-1}\sigma \check{Z}(K_0)\big)\exp\frac{f}{2}\omega_1\big)\Big) &
\end{xalignat*}

We have the following fact:
\begin{xalignat*}{2}
W_{\mathfrak{g}}\big(\hat{\Omega}^{-1}\sigma \check{Z}(K_0)\big)(\lambda)=&W_{\mathfrak{g}}\big(\hat{\Omega}^{-1}\sigma \hat{Z}(K_0)\big)(\lambda)W_{\mathfrak{g}}\big(\hat{\Omega}^{-1}\sigma(\chi (\Omega))\big)(\lambda) & \text{($\hat{\Omega}^{-1}\sigma$ : an alg. iso.)} \\
=&D(j_{\mathfrak{g}} ^{\frac{1}{2}})^{-1}W_{\mathfrak{g}}\sigma \hat{Z}(K_0)(\lambda) D(j_{\mathfrak{g}} ^{\frac{1}{2}})^{-1}W_{\mathfrak{g}}\sigma \chi \Omega(\lambda) & \text{(Theorem \ref{CM})} \\
=&W_{\mathfrak{g}}\sigma \hat{Z}(K_0)(\lambda-\rho)W_{\mathfrak{g}}\sigma \hat{Z}(\bigcirc)(\lambda-\rho) & \text{(by Remark \ref{Du})}\\
=& \frac{\prod _{\alpha >0}(\rho,\alpha)^2}{\prod _{\alpha >0}(\lambda,\alpha)^2}Q_\mathfrak{g}(K_0)\big|_{q=e^h} \big(\lambda-\rho) & \text{(section \ref{s-Q})} & \\*
=&\frac{\prod _{\alpha >0}(\rho,\alpha)^2}{\prod _{\alpha >0}(\lambda,\alpha)^2}\sum_{\beta \in Y ,n\in \mathbb{Z}_+ \atop {2|\Phi_+|\le j\le n+2|\Phi_+|}} c_{\beta,j,n}\beta^j(\lambda)h^n &
\end{xalignat*}

The following equation is known, where $W$ is the Weyl group.
$$\prod _{\alpha >0}(q^{\frac{\alpha}{2}}-q^{-\frac{\alpha}{2}})=\sum_{w \in W}\mathrm{sign}(w)q^{w(\rho)}$$
Then we have
\begin{xalignat*}{2}
\bigg(\prod _{\alpha >0}(q^{\frac{\alpha}{2}}-q^{-\frac{\alpha}{2}})\bigg)^2
&=\bigg( \sum_{w \in W}\mathrm{sign}(w)q^{w(\rho)} \bigg)^2 & \\
&=\sum_{w,w' \in W}\mathrm{sign}(ww')q^{w(\rho)+w'(\rho)} &
\end{xalignat*}

To a vector $x$ we associate the directional derivation $\partial _x$, whose action is given by
$\partial _x \lambda=(\lambda,x)$.
Note 
$$\sum _{k=1} ^{\mathrm{dim}\mathfrak{g}}(x_k ,\lambda)^2=(\lambda,\lambda)$$
where $\{x_k\}$ is an orthnormal basis with respect to the form $(\ ,\ )$.

\begin{xalignat*}{2}
&\exp(-\frac{h}{2f}\sum_k\partial _{x_k}^2) \bigg(\sum_{w,w' \in W}\mathrm{sign}(ww')q^{w(\rho)+w'(\rho)}\bigg) \bigg|_{\alpha\to h^{-1}\alpha}\bigg|_{x=0}& \\
=&\sum_{w,w' \in W}\mathrm{sign}(ww')\sum_{j=0} ^{\infty } \big(-\frac{h}{2f}\big)^j\frac{1}{j!}(\sum_k\partial _{x_k}^2)^{j}\bigg(\sum_{i=0} ^{\infty }\frac{1}{i!}\big(w(\rho)+w'(\rho)\big)^i \bigg) \bigg|_{x=0}& \\
=&\sum_{w,w' \in W}\mathrm{sign}(ww')\sum_{j=0} ^{\infty }\big(-\frac{h}{2f}\big)^j\frac{1}{j!}\frac{1}{(2j)!}(2j)! \big(w(\rho)+w'(\rho),w(\rho)+w'(\rho)\big)^j & \\
=&\sum_{w,w' \in W}\mathrm{sign}(ww')\sum_{j=0} ^{\infty }\big(-\frac{h}{2f}\big)^j\frac{1}{j!}\big(2|\rho|^2+2(w(\rho),w'(\rho))\big)^j & \\
=&\sum_{w,w' \in W}\mathrm{sign}(ww')\sum_{j=0} ^{\infty }\big(-\frac{h}{f}\big)^j\frac{1}{j!}\big(|\rho|^2+(w(\rho),w'(\rho))\big)^j & \\
=&q^{-\frac{|\rho|^2}{f}} \sum_{w,w' \in W}\mathrm{sign}(ww')q^{(w(\rho),w'(\rho))}  & \\
=&|W|q^{-\frac{|\rho|^2}{f}} \sum_{w \in W}\mathrm{sign}(w)q^{(\rho,w(\rho))} & \\
=&|W|q^{-\frac{|\rho|^2}{f}}\prod _{\alpha >0}(q^{\frac{(\rho,\alpha)}{2f}}-q^{-\frac{(\rho,\alpha)}{2f}}) & \\
=&|W|\prod _{\alpha >0}(1-q^{-\frac{(\rho,\alpha)}{f}}) & \\
=&|W|\bigg(\big(-\frac{h}{f}\big)^{|\Phi_+|}\prod _{\alpha >0}(\rho,\alpha)\bigg) \times \big(1+ \cdot\cdot\cdot \big) &
\end{xalignat*}

From the above calculation, we know that $|W|\prod _{\alpha >0}(\rho,\alpha)\big(-\frac{h}{f}\big)^{|\Phi_+|}$ is the result of action on $\prod _{\alpha >0}\alpha ^2$ by $\exp(-\frac{h}{2f}\sum_k\partial _{x_k}^2)$.

\begin{xalignat*}{2}
&\hat{W}_{\mathfrak{g}}\Big(\int^{FG} \big(\hat{\Omega}^{-1}\sigma \check{Z}(K_0)\big)\exp\frac{f}{2}\omega_1\big)\Big) & \\
=&\hat{W}_{\mathfrak{g}}\big(\bigg\langle\exp(-\frac{\omega_1}{2f}),\hat{\Omega}^{-1}\sigma \check{Z}(K_0) \bigg\rangle \big)  & \\
=&\exp(-\frac{h}{2f}\sum_k\partial _{x_k}^2) \big(\prod _{\alpha >0}(\rho,\alpha)^2\sum c_{\beta,i,n}\frac{1}{\prod _{\alpha >0}(h^{-1}\alpha )^2 }h^n (h^{-1}\beta)^i ) \big) \bigg|_{x=0} & \text{}\\
=&h^{2|\Phi_+|}\prod _{\alpha >0}(\rho,\alpha)^2\sum_j \big(-\frac{h}{2f}\big)^j\frac{1}{j!}(\sum_k\partial _{x_k}^2)^{j}\bigg(\sum c_{\beta,i,n}\frac{1}{\prod _{\alpha >0}\alpha ^2 }h^n (h^{-1}\beta)^i ) \bigg) \bigg|_{x=0} & \\
=&h^{2|\Phi_+|}\prod _{\alpha >0}(\rho,\alpha)^2\sum c_{\beta,2j,n}\frac{1}{|W|\prod _{\alpha >0}(\rho,\alpha)} \big(-\frac{f}{h} \big)^{|\Phi_+|}\big(-\frac{|\beta|^2}{f}\big)^j h^{n-j}\frac{1}{2^j \  j!} (2j)! & \\*
=&(-fh)^{|\Phi_+|}\frac{\prod _{\alpha >0}(\rho,\alpha)}{|W|} \sum c_{\beta,2j,n}\big(-\frac{|\beta|^2}{f}\big)^j (2j-1)!! h^{n-j} & 
\end{xalignat*}
In the second equality, we identify the contraction map (see section \ref{WS}) with the directional derivation.

By (\ref{Omega}) and (\ref{Theta}), 
\begin{xalignat*}{2}
&\hat{W}_{\mathfrak{g}}(\langle\Omega,\Omega\rangle)\exp \big( \frac{3\mathrm{sign}(f)-f}{48}\hat{W}_{\mathfrak{g}}(\theta) \big) & \\
=&\bigg( \prod _{\alpha >0} \frac{\text{sinh} \frac{(\rho,\alpha)}{2}h} {\frac{(\rho ,\alpha )}{2}h} \bigg)(e^h)^{\frac{3\mathrm{sign}(f)-f}{2}|\rho|^2} & \\*
=&\frac{1}{\prod_{\alpha>0}(\rho,\alpha)}\Big(-\frac{\mathrm{sign}(f)}{h}\Big)^{|\Phi_+|}q^{\frac{\mathrm{sign}(f)-f}{2}|\rho|^2}\prod_{\alpha>0} (1-q^{\mathrm{sign}(f)(\rho,\alpha)}) & 
\end{xalignat*}

Recall the following fact:
$$|H_1(M_K;\mathbb{Z})|=f\mathrm{sign}(f)$$

Thus we have
$$\hat{W}_{\mathfrak{g}}\big( \hat{Z}^{\mathrm{LMO}}(M_K) \big)=|H_1(M_K;\mathbb{Z})|^{|\Phi_+|}\tau^{PG}(M_K)$$
In a general case, by theorem \ref{37}, we can assume $L$ is an algebraiclly split link. In this case, the proof is similar to above proved-case, using $\lq\lq$partial differential".
\end{proof}


\begin{thebibliography}{9999}
 \baselineskip18pt               
\bibitem{BN} D. Bar-Natan,
{\it On the Vassiliev invariants }, Topology, {\bf 34} (1995), 423--472.

\bibitem{Wheel} D. Bar-Natan, S. Garoufalidis, L. Rozansky and D.P. Thurston,
{\it Wheels, wheeling, and the Kontsevich integral of the unknot }, Israel J.Math. {\bf 119} (2000), 217--238.

\bibitem{BGRT} D. Bar-Natan, S. Garoufalidis, L. Rozansky and D.P. Thurston,
{\it The $\r{A}$rhus integral of rational homology 3-spheres. I. A highly non trivial flat connection on $S\sp 3$ }, Selecta Math. (N.S.) {\bf 8} (2002), no. 3, 315--339.

\bibitem{BGRT2} D. Bar-Natan, S. Garoufalidis, L. Rozansky and D.P. Thurston,
{\it The $\r{A}$rhus integral of rational homology 3-spheres. II. Invariance and universality }, Selecta Math. (N.S.) {\bf 8} (2002), no. 3, 341--371.

\bibitem{BGRT3} D. Bar-Natan, S. Garoufalidis, L. Rozansky and D.P. Thurston,
{\it The $\r{A}$rhus integral of rational homology 3-spheres III: The Relation with the Le-Murakami-Ohtsuki Invariant }, Selecta Math. (N.S.) {\bf 10} (2004), no. 3, 305--324.

\bibitem{BL} D. Bar-Natan and R. Lawrence,
{\it A rational surgery formula for the LMO invariant }, Israel J. Math. {\bf 140} (2004), 29--60.

\bibitem{BGV} N. Berline, E. Getzler and M. Vergne ,
{\it Heat kernels and Dirac operators },
Grundlehen \ der \ mathematischen \ wissenschaften {\bf 298} , Springer-Verlag \ Berlin \ Heidelberg \ 1992.

\bibitem{G0} S. Garoufalidis,
{\it Rationality:From Lie algebras to Lie groups }, preprint.

\bibitem{G1} S. Garoufalidis,
{\it Beads:From Lie algebras to Lie groups }, preprint.

\bibitem{Ha1} K. Habiro,
{\it Cyclotomic completions of polynomial rings }, Publ. Res. Inst. Math. Sci, {\bf 40} (2004), no. 4, 1127--1146.

\bibitem{Ha2} K. Habiro,
{\it On the quantum $\mathfrak{sl_2}$ invariants of knots and integral homology spheres }, Geometry and Topology Monographs, Vol. 4 (2002), 55-68.

\bibitem{TK} T. Kuriya,
{\it The LMO invariant and Guadagnini-Pilo's conjecture for lens spaces }, preprint.

\bibitem{LMO} T.T.Q. Le , J.Murakami and T.Ohtsuki,
{\it On a universal perturbative invariant of 3-manifolds },
Topology, {\bf 37} (1998), 539--574.

\bibitem{Le2} T.T.Q. Le ,
{\it Quantum invariants of 3-manifolds: integrality, splitting, and perturbative expansion }, Topology Appl. {\bf 127} (2003), no. 1-2, 125--152.

\bibitem{Oh} T. Ohtsuki ,
{\it Quantum invariants ,A study of knots, 3-manifolds, and their sets },
Series \ on \ Knot \ and \ Everything {\bf 29} , World \ Scientific \ Publishing \ Co., Inc \ 2001.

\bibitem{Oh2} T. Ohtsuki ,
{\it The perturbative $SO$(3) inavriant of rational homology 3-spheres recovers from the universal perturbative invariant },
Topology {\bf 39} , (2000) \ 1103--1135.

\bibitem{OH3} Edited by T. Ohtsuki ,
{\it Problems on invariants of knots and 3--manifolds  },
Geom. Topol. Publ., Coventry, 2002.

\bibitem{RT} N. Yu. Reshetikhin, V. Turaev,
{\it Invariants of 3-manifolds via link polynomials and quantum groups},
Invent. Math., {\bf 103} (1991), 547--598.

\bibitem{Ro} L. Rozansky,
{\it A universal U(1)-RCC invariant of links and rationality conjecture }, preprint, 2002 \ math.GT/0201139.

\bibitem{Ro2} L. Rozansky,
{\it A rationality conjecture about Kontsevich integral of knots and its implications to the structure of the colored Jones polynomial }, Topology Appl. {\bf 127} (2003), no. 1-2, 47--76.

\bibitem{Wi} S. Willerton,
{\it The Kontsevich Integral and structures on the space of diagrams }, Knots in Hellas '98 (Delphi), 530--546, Ser. Knots Everything, 24, World Sci. Publishing, River Edge, NJ, 2000.

\bibitem{V} A. Vaintrob,
{\it Universal weight systems and the Melvin-Morton expansion of the colored Jones knot invariant }, Algebraic geometry, 5. J. Math. Sci. {\bf 82} (1996), no. 1, 3240--3254.

\end{thebibliography}
\end{document}